\documentclass[numreferences]{kluwer}

\usepackage{amsthm}
\usepackage{amsmath}
\usepackage{amssymb}
\usepackage{graphicx}
\usepackage{epsfig}
\usepackage{graphicx}

\newtheorem{theorem}{Theorem}

\newtheorem{definition}[theorem]{Definition}

\newtheorem{corollary}[theorem]{Corollary}
\newtheorem{remark}[theorem]{Remark}
\newtheorem{example}{Example}
\numberwithin{equation}{section}

\numberwithin{theorem}{section}

\begin{document}
\begin{opening}

\title{Parameter estimation in diagonalizable bilinear \\stochastic
parabolic equations}

\author{Igor \surname{Cialenco} }
\institute{Department of Applied Mathematics,
Illinois Institute of Technology\\
10 West 32nd Str, Bld E1, Room 208,
Chicago, IL 60616, USA\\
\email{igor@math.iit.edu, http://math.iit.edu/$\sim$igor} }

\author{Sergey V. \surname{Lototsky}  \thanks{SVL acknowledges support
from the NSF CAREER award DMS-0237724.} }
\institute{Department of Mathematics, University of Southern California\\
3620 S. Vermont Avenue, KAP 108,
Los Angeles, CA 90089, USA \\
\email{lototsky@math.usc.edu, http://math.usc.edu/$\sim$lototsky} }

\runningtitle{Estimation in bilinear equations}
\runningauthor{Ig. Cialenco, S. V. Lototsky}
\classification{AMS 2000}{Primary 62F12; Secondary  60H15}

\begin{abstract}
A parameter estimation problem is considered for a
 stochastic parabolic equation with multiplicative noise under the
 assumption that the equation can be reduced to an infinite
 system of uncoupled diffusion processes. From the point of view
 of classical statistics, this problem turns out to be singular not
 only for the original infinite-dimensional system but also for
 most finite-dimensional projections. This singularity can be
 exploited to improve the rate of convergence of traditional estimators
 as well as to construct completely new closed-form exact estimator.
\end{abstract}

\end{opening}
%\today

\section{Introduction}

In the classical statistical estimation problem, the starting point is a family
$\mathbf{P}^{\theta}$ of probability measures depending on the parameter
$\theta$ belonging to some subset $\Theta$ of a finite-dimensional
Euclidean space. Each $\mathbf{P}^{\theta}$ is the distribution of a random
element. It is assumed that a  realization of one random element
 corresponding to one value $\theta=\theta_0$ of the parameter is
 observed, and the objective is to estimate the values of this
 parameter from the observations.

 The intuition is to select the value $\theta$ corresponding to
 the random element that is {\em most likely} to produce the observations.
 A rigorous mathematical implementation of this idea leads to the notion
 of the regular statistical model \cite{IbragimovKhasminskiiBook1981}:
 the statistical model (or estimation problem)
 $\mathbf{P}^{\theta},\ \theta\in\Theta$, is
 called regular, if the following two conditions are satisfied:
 \begin{itemize}
 \item there exists a probability measure $\mathbf{Q}$ such that
 all measures $\mathbf{P}^{\theta}$ are absolutely continuous
 with respect to $\mathbf{Q}$;
 \item the density $d\mathbf{P}^\theta/d\mathbf{Q}$, called the
 likelihood ratio, has a special property,
 called local asymptotic normality.
 \end{itemize}
 If at least one of the above conditions is violated, the problem is
 called singular.

 In regular models, the estimator $\widehat{\theta}$
  of the unknown parameter is constructed by
 maximizing the likelihood ratio and is called the maximum likelihood estimator
 (MLE). Since, as a rule, $\widehat{\theta}\not=\theta_0$,
 the consistency of the estimator is studied, that is, the convergence
 of $\widehat{\theta}$ to $\theta_0$ as more and more information becomes
 available. In all known regular statistical problems, the amount of information
 can be increased in one of two ways: (a) Increasing the sample size,
 for example, the observation time interval (large sample asymptotic);
(b) reducing the amplitude of noise (small noise asymptotic). The
asymptotic behavior of $\widehat{\theta}$ in both cases is well-studied.
It is also known that many other estimators in regular
models are asymptotically equivalent to the MLE.

While all regular models are in a sense the same, each singular model is
different. Sometimes, it is possible to approximate
 a singular model  with a sequence of regular
models. For each regular model, an MLE is constructed, and then in the limit
one can often get the true value of the parameter while both the sample size
and the noise amplitude are fixed. Some singular models cannot be approximated
by a sequence of regular models and admit estimators that have nothing to do
with the MLE \cite{KhasminskiiKrylovMoshniuk}.
 In this paper, Section \ref{subsection-exact}, we introduce
 a completely new type of such estimators for a large class of
 singular models.

Infinite-dimensional stochastic evolution equations, that is,
stochastic evolution equations in infinite-dimensional spaces,
are a rich source of statistical problems, both regular and singular.
A typical example is the It\^{o} equation
\begin{equation}\label{eq1}
\begin{cases}
  du(t) + (\mathcal{A}_0 + \theta \mathcal{A}_1) u(t)dt =
  f(t)dt+ \sum\limits_{j\geq 1}(\mathcal{M}_ju(t)+g_j(t))dW_j(t), \\
  u(0) = u_0,
\end{cases}
\end{equation}
where $t\in[0,T], \ \mathcal{A}_0, \ \mathcal{A}_1, \ \mathcal{M}_j$
are linear operators,
$f, \ g_j$ are adapted processes,
$W_j$ are independent Wiener processes,
 and $\theta$ is the unknown parameter belonging to an open
subset of the real line. The underlying assumption is that the
solution $u$ exists, is unique, and
 can be observed as an infinite-dimensional object for
all $t\in [0,T]$. Depending on the operators in the equation, the
estimation model can be regular, a singular limit of regular problems, or
completely singular.

If $\mathcal{A}_0, \mathcal{A}_1, \ \mathcal{M}_j$ are  partial differential or
pseudo-differential operators, \eqref{eq1} becomes a  stochastic partial
differential equation (SPDE), which is becoming increasingly popular
for  modelling various phenomena in
 fluid mechanics \cite{SerranoUnny1990},
 oceanography \cite{Piterbarg2001},
temperature  anomalies \cite{Frankignoul2000,Piterbarg2006},
finance \cite{Cont2005, GallPap, Goldstein00},  and other domains.
Various estimation problems for different types of SPDEs have been
investigated by many authors:
\cite[etc.]{Aihara92, BagchiBorkar, HuebnerLototskyRozovskii97,
HubnerRozovskiiKhasminskii, HuebnerRozovskii, IbragimovKhasminskii1998,
IbragimovKhasminskii1999, IbragimovKhasminskii2000, Lototsky2003,
LototskyRozovskii99, LototskyRozovskii2000}.

Depending on the stochastic part,  \eqref{eq1} is classified  as follows:
\begin{itemize}
\item equation with {\bf additive noise}, if $\mathcal{M}_j=0$ for all $j$;
 \item equation  with {\bf multiplicative noise} (or bilinear equation) otherwise.
 \end{itemize}
 Depending on the operators, \eqref{eq1} is classified as follows:
 \begin{itemize}
 \item {\bf Diagonalizable equation}, if the operators $\mathcal{A}_0$,
 $\mathcal{A}_1$, and $\mathcal{M}_j$, $j\geq 1$, have a common
 system of eigenfunctions $h_k,\ k\geq 1$, and this system is an
 orthonormal basis in a suitable Hilbert space.
 \item {\bf Non-diagonalizable equation}  otherwise.
 \end{itemize}
A diagonalizable equation is reduced to an infinite system of
{\em uncoupled} one-dimensional diffusion processes; these processes are
the Fourier coefficients of the solution in the basis $h_k$. As a result,
 while somewhat restrictive as a modelling tool, diagonalizable equations
 are an extremely convenient object to  study  estimation problems and often
 provide the benchmark results that carry over to more general equations.

The parameter estimation problem for a diagonalizable equation
\eqref{eq1} with additive space-time white noise (that is,
$g_j=h_j$ and $\mathcal{M}_j=0$ for all $j$)
 was studied for the first time by
 Huebner, Khasminskii, and Rozovskii \cite{HubnerRozovskiiKhasminskii},
and further investigated in  \cite{HuebnerLototskyRozovskii97,
HubnerRozovskiiKhasminskii, HuebnerRozovskii, PiterbargRozovskii1997}.
The main feature of this problem is that every $N$-dimensional
projection of the equation leads to a regular statistical problem, but
the problem can become singular in the limit $N\to \infty$
(a singular limit of regular problems); when
this happens, the dimension $N$ of the projection becomes a natural
asymptotic parameter of the problem. Once the diagonalizable
model is well-understood, extensions to more general equations can
be considered (\cite{LototskyRozovskii99, LototskyRozovskii2000}).

This paper is the first attempt to investigate the estimation problem
for infinite-dimen\-sional bilinear equations. Such models are often
completely singular, that is, cannot be represented as a limit
of regular models.
We consider the more tractable situation of diagonalizable equations.
In Section \ref{Section-Settings} we provide the necessary background
on stochastic evolution equations, with emphasis on diagonalizable
bilinear equations. The maximum likelihood estimator (MLE)
and its modifications for  diagonalizable
bilinear equations are studied in Section \ref{section-estimates}.
  We give sufficient conditions  on operators
$\mathcal{A}_0, \mathcal{A}_1, \mathcal{M}$ that ensure  consistency and asymptotic
 normality of the MLE. We also demonstrate that the MLE in this model is
 not always the best estimator, which, for a singular model is
 not at all surprising. Section \ref{subsection-exact} emphasizes the point
 even more by introducing  a {\tt closed-form exact estimator}. Due to the
specific structure of stochastic term, for a large class of
infinite-dimensional systems with finite-dimensional noise, one can get
the {\em exact} value of the unknown parameter after a finite number
of arithmetic manipulations with the observations. The very existence of
such estimators in these models is rather remarkable and has no analogue
in classical statistics.

As an illustration,  let $\theta$ be a positive number, $W$ a standard Wiener
process,  and consider the It\^{o} equation
\begin{equation}
\label{intE}
du(t,x) - \theta u_{xx}(t,x)dt = u(t,x)dW(t),\ t>0,\ x\in (0,\pi),
\end{equation}
with zero boundary conditions. If
$h_k(x)=\sqrt{2/\pi}\, \sin(kx)$, $k\geq 1,$ and
$$
u_k(t)=\int_0^{\pi} u(t,x)h_k(x)dx \ ,
$$
then
\begin{equation}
\label{intrFS}
u(t,x)=\sum_{k\geq 1}u_k(t)h_k(x)
\end{equation}
 and each $u_k$ is a geometric Brownian motion:
$$
u_k(t)=u_k(0)-\int_0^tk^2u_k(s)ds+\int_0^tu_k(s)dW(s).
$$
We assume that $u_k(0)\not=0$ for all $k\geq 1$.
In Sections \ref{section-estimates}
and \ref{subsection-exact} we establish the following result.

\begin{theorem}
If the solution of equation \eqref{intE} is observed in the form
\eqref{intrFS}, then  the parameter $\theta$ can be computed in
each of the following ways:
\begin{enumerate}
\item[(E1)] $\displaystyle
\theta=\lim_{T\to \infty} \left(\frac{1}{k^2T}\ln\frac{u_k(0)}{u_k(T)}-
\frac{1}{2k^2}\right)$ for every $k\geq 1$;
\item[(E2)] $\displaystyle
\theta=\lim_{k\to \infty} \frac{1}{k^2T}\ln\frac{u_k(0)}{u_k(T)}$
for every $T>0$;
\item[(E3)] $\displaystyle
 \theta=\frac{1}{T(k^2-n^2)}\ln\frac{u_n(T)u_k(0)}{u_k(T)u_n(0)}$
for every $T>0$ and $n\not= k$.
\end{enumerate}
\end{theorem}

Both (E1) and (E2) are essentially the same
maximum likelihood estimator, but the infinite-dimensional
nature of the equation makes it possible to study this estimator in two
different asymptotic regimes. (E3) is a closed-form exact estimator.
While it is most likely to be the best choice for this particular problem,
 we show in Section \ref{subsection-exact} that computational
complexity of closed-form exact estimators can dramatically increase
with the number of Wiener processes driving the equation, while
the complexity of the MLE is almost unaffected by this number. The result is
another unexpected feature of closed-form exact estimators: ever though
they produce the exact value of the parameter, they are
not always the best choice computationally.

\section{Stochastic Parabolic Equations}
\label{Section-Settings}
In this section we introduce the diagonalizable stochastic parabolic equation
depending on a parameter and study the main properties of the
solution.

Let $\mathbf{H}$ be a separable Hilbert space with the inner product $(\cdot,
\cdot)_0$ and the corresponding norm $\|\cdot\|_0$. Let $\Lambda$ be  a
densely-defined
linear operator on $\mathbf{H}$ with the following property:
there exists a positive number $c$ such that
 $\| \Lambda u\|_0\geq c\|u\|_0$ for
every $u$ from the domain of $\Lambda$.
Then the operator powers  $\Lambda^\gamma, \ \gamma\in\mathbb R,$ are
 well defined and generate the spaces $\mathbf{H}^\gamma$:
 for $\gamma>0$, $\mathbf{H}^\gamma$ is the domain of
  $\Lambda^\gamma$; $\mathbf{H}^0=\mathbf{H}$;
 for $\gamma <0$,  $\mathbf{H}^\gamma$ is the completion
 of $\mathbf{H}$ with respect to the norm
 $\| \cdot \|_\gamma := \|\Lambda \cdot\|_0$ (see for instance
 Krein at al. \cite{KreinPetuninSemenov}).
 By construction, the collection of spaces
 $\{ \mathbf{H}^\gamma,\ \gamma\in \mathbb R\}$
  has the following properties:
 \begin{itemize}
 \item $\Lambda^{\gamma}(\mathbf{H}^r) = \mathbf{H}^{r-\gamma}$ for every
$\gamma,r\in\mathbb{R}$;
 \item For $\gamma_1<\gamma_2$ the space $\mathbf{H}^{\gamma_2}$ is densely and
 continuously embedded into $\mathbf{H}^{\gamma_1}$: $\mathbf{H}^{\gamma_2}
 \subset \mathbf{H}^{\gamma_1}$ and there exists a positive
 number $c_{12}$ such that  $\|u\|_{\gamma_1}\leq c_{12} \|u\|_{\gamma_2}$
 for all $u\in \mathbf{H}^{\gamma_2}$ ;
\item  for every
$\gamma\in\mathbb R$ and $m>0$, the space $\mathbf{H}^{\gamma -m}$ is the
dual of $\mathbf{H}^{\gamma+m}$ relative to the inner product in
$\mathbf{H}^{\gamma}$, with duality $\langle\cdot,\cdot\rangle_{\gamma,m}$ given by
$$
\langle u_1, u_2 \rangle_{\gamma,m} = (\Lambda^{\gamma -m}u_1,
\Lambda^{\gamma+m}u_2)_0, \ {\rm where\ }
u_1\in\mathbf{H}^{\gamma-m},\ u_2\in\mathbf{H}^{\gamma +m}.
$$
\end{itemize}

Let $(\Omega, \mathcal{F}, \{\mathcal{F}_t\}, \mathbb P)$ be a stochastic basis
with the usual assumptions, and let $\{W_j,\ j\geq 1\}$ be a
collection of independent standard Brownian motions on this basis.
 Consider the following It\^{o} equation
\begin{equation}
\label{eq2}
\begin{cases}
  du(t) + (\mathcal{A}_0 + \theta \mathcal{A}_1) u(t)dt
  = f(t)dt +
  \sum\limits_{j\geq 1}(\mathcal{M}_ju(t)+g_k(t))dW_j(t), \ 0<t\leq T,  \\
  u(0) = u_0\,
\end{cases}
\end{equation}
where $\mathcal{A}_0, \, \mathcal{A}_1, \, \mathcal{M}_j$ are linear operators,
$f$ and $g_k$ are  adapted process,
 and $\theta$ is a scalar parameter  belonging to
an open set $\Theta\subset \mathbb{R}$.

\begin{definition}
\label{def000}
{\ \ }\\
(a) Equation  \eqref{eq2} is called an {\tt equation with additive noise} if
$\mathcal{M}_j=0$ for all $j\geq 1$. Otherwise,
\eqref{eq2} is called an {\tt equation with multiplicative noise} (also known as
a {\tt bilinear equation}). \\
(b) Equation  \eqref{eq2} is called  {\tt diagonalizable} if
the  operators $\mathcal{A}_0,\ \mathcal{A}_1, \mathcal{M}_j,\ j\geq 1$,
have a common system of eigenfunctions $\{h_k,\ k\geq 1\}$ such that
$\{h_k,\ k\geq 1\}$ is an orthonormal basis in $\mathbf{H}$ and
each $h_k$ belongs to every $\mathbf{H}^{\gamma}$. \\
(c) Equation  \eqref{eq2} is called {\tt parabolic  in the triple}
$(\mathbf{H}^{\gamma+m}, \mathbf{H}^{\gamma},\mathbf{H}^{\gamma-m})$ if
\begin{itemize}
\item[$\circ$]
the operator $\mathcal{A}_0+\theta\mathcal{A}_1$ is uniformly bounded from
$\mathbf{H}^{\gamma+m}$ to $\mathbf{H}^{\gamma-m}$ for
$\theta \in \Theta:$ there exists a positive real number $C_1$
such that
\begin{equation}
\label{contA}
\|(\mathcal{A}_0+\theta\mathcal{A}_1)v\|_{\gamma-m}\leq C_1\|v\|_{\gamma+m}
\end{equation}
for all $\theta\in \Theta$, $v\in \mathbf{H}^{\gamma+m}$;
\item[$\circ$]
 There exists a positive number $\delta$ and a real number $C$
such that, for every $v\in \mathbf{H}^{\gamma+m}$, $\theta\in \Theta$,
\begin{equation}
\label{parab}
-2\langle (\mathcal{A}_0+\theta\mathcal{A}_1)v,v\rangle_{\gamma,m}
+\sum_{j\geq 1} \|\mathcal{M}_jv\|_{\gamma}^2+\delta\|v\|_{\gamma+m}^2\leq
C\|v\|_{\gamma}^2.
\end{equation}
\end{itemize}
\end{definition}

\begin{remark}
\label{rm00}
{\rm
(a) Note that \eqref{contA} and \eqref{parab} imply uniform continuity of
the family of operators $\mathcal{M}_j$, $j\geq 1$ from
$\mathbf{H}^{\gamma+m} $ to $\mathbf{H}^{\gamma}$; in fact,
$$
\sum_{j\geq 1} \|\mathcal{M}_jv\|_{\gamma}^2 \leq
2C_1\|v\|_{\gamma+m}^2+C\|v\|_{\gamma}^2.
$$
(b) If equation \eqref{eq2} is parabolic, then condition \eqref{parab}
implies that
$$
\langle (2\mathcal{A}_0+2\theta\mathcal{A}_1+CI)v,v\rangle_{\gamma,m}
\geq \delta\|v\|_{\gamma+m}^2,
$$
where $I$ is the identity operator.
The Cauchy-Schwartz inequality and the continuous embedding of
$\mathbf{H}^{\gamma+m}$ into $\mathbf{H}^{\gamma}$ then imply
$$
\|(2\mathcal{A}_0+2\theta\mathcal{A}_1+CI)v\|_{\gamma}\geq \delta_1\|v\|_{\gamma}
$$
for some $\delta_1>0$ uniformly in $\theta\in \Theta$.
As a result, we can take $\Lambda=
(2\mathcal{A}_0+2\theta^*\mathcal{A}_1+CI)^{1/(2m)}$ for some fixed
$\theta^*\in\Theta$.
}
\end{remark}

{\em From now on, if equation \eqref{eq2} is parabolic and diagonalizable,
we will assume that the operator $\Lambda$ has the same eigenfunctions as the
operators $\mathcal{A}_0,\ \mathcal{A}_1,\ \mathcal{M}_j$;
by Remark \ref{rm00}, this leads to no loss of generality.}

\begin{example}
\label{ex:main}
{\rm
(a) For $0<t\leq T$ and $x\in (0,1)$, consider the equation
\begin{equation}
\label{ex00}
du(t,x)-\theta\,u_{xx}(t,x)dt= u_x(t,x)dW(t)
\end{equation}
with periodic boundary conditions; $u_x=\partial u/\partial x$.
Then $\mathbf{H}^{\gamma}$ is the Sobolev space on the
unit circle (see, for example, Shubin \cite[Section I.7]{Shubin}) and
$\Lambda=\sqrt{I-\boldsymbol{\Delta}}$, where $\boldsymbol{\Delta}$ is
the Laplace operator on $(0,1)$ with periodic boundary conditions.
Direct computations show that equation \eqref{ex00} is diagonalizable;
it is parabolic if and only if $2\theta>1$.

(b) Let $G$ be a smooth bounded domain in
$\mathbb{R}^d$. Let $\boldsymbol{\Delta}$ be the Laplace operator
on $G$ with zero boundary conditions. It is known  {\rm (for example, from
Shubin \cite{Shubin})}, that
\begin{enumerate}
\item the eigenfunctions $\{h_k,\ k\geq 1\}$
 of $\boldsymbol{\Delta}$
are smooth in $G$ and form an orthonormal basis in $L_2(G)$;
\item the corresponding eigenvalues $\sigma_k,\ k\geq 1$,
can be arranged so that $0<-\sigma_1\leq -\sigma_2\leq \ldots$, and there
exists a number $c>0$ such that $|\sigma_k|\sim ck^{2/d}$, that is,
$$
\lim_{k\to \infty} |\sigma_k|k^{-2/d}=c.
$$
\end{enumerate}

We take $\mathbf{H}=L_2(G)$,
$\Lambda=\sqrt{I-\boldsymbol{\Delta}}$, where $I$ is the identity operator.
 Then $\|\Lambda u\|_0\geq \sqrt{1-\sigma_1}
\|u\|_0$ and the operator $\Lambda$ generates the Hilbert spaces
$\mathbf{H}^{\gamma}$, and,  for every $\gamma\in \mathbb{R}$, the space
$\mathbf{H}^{\gamma}$ is the closure of the set of smooth compactly supported
function on $G$ with respect to the norm
$$
\left(\sum_{k\geq 1} (1+k^2)^{\gamma}|\varphi_k|^2\right)^{1/2},\
{\rm \ where\ } \varphi_k=\int_G \varphi(x)h_k(x)dx,
$$
which is an equivalent norm in $\mathbf{H}^{\gamma}$.
  Let $\theta$ and $\sigma$ be real numbers.
Then the stochastic equation
\begin{equation}
\label{eq:exmain}
du-\theta\boldsymbol{\Delta}udt=\Lambda u\,dW
\end{equation}
is
\begin{itemize}
\item always diagonalizable;
\item parabolic in $(\mathbf{H}^{\gamma+1},\mathbf{H}^{\gamma},
\mathbf{H}^{\gamma-1})$ for every $\gamma\in \mathbb{R}$
  if and only  if
$2\theta>1$.
\end{itemize}
Indeed, we have $\mathcal{A}_0=0$, $\mathcal{A}_1=-\boldsymbol{\Delta}$,
$\mathcal{M}_1=\Lambda$, $\mathcal{M}_j=0$, $j\geq 2$, and
$$
-2\theta\langle \mathcal{A}_1v,v\rangle_{\gamma,1}=-2\theta\|v\|_{\gamma+1}^2+
2\theta\|u\|_{\gamma}^2,
$$
and so \eqref{parab} holds with $\delta=2\theta-1$ and $C=2\theta$.
}
\end{example}

\begin{remark}
{\rm Taking in \eqref{eq2} $\mathbf{H}=L_2(G)$, where $G$ is a
smooth bounded domain in $\mathbb{R}^d$, and
 $\mathcal{A}_0=-\boldsymbol{\Delta}$, $\mathcal{A}_1=I$,
 $\mathcal{M}_ju=h_k(x)u(x)$, $g_k=h_k(x)g(t,x)$,
we get a bilinear  equation driven by space-time white noise.
Direct analysis shows that this equation is not diagonalizable. Moreover, the
equation is parabolic if and only if $d=1$, that is, when $G$ is an interval;
for details, see the lecture notes by Walsh \cite{Walsh}.
}
\end{remark}

For a diagonalizable equation, the parabolicity condition
\eqref{parab} can be expressed in terms of the eigenvalues of the
operators in the equation.

\begin{theorem}
\label{th0}
Assume that equation \eqref{eq2} is diagonalizable, and
$$
\mathcal{A}_0h_k=\rho_kh_k,\
\mathcal{A}_1h_k=\nu_kh_k,\
\mathcal{M}_jh_k=\mu_{jk}h_k.
$$
With no loss of generality {\rm (see Remark \ref{rm00})},
 we also assume that
 $$
 \Lambda h_k=\lambda_kh_k.
 $$
 Then equation \eqref{eq2} is parabolic in the triple
$(\mathbf{H}^{\gamma+m}, \mathbf{H}^{\gamma},\mathbf{H}^{\gamma-m})$
if and only if there exist  positive real numbers $\delta, C_1$ and
a real number $C_2$ such that, for all $k\geq 1$ and $\theta\in \Theta$,
\begin{align}
&\lambda_k^{-2m}|\rho_k+\theta\nu_k|  \leq C_1; \label{eig1} \\
&-2(\rho_k+\theta\nu_k)+ \sum_{j\geq 1} |\mu_{jk}|^2 +
\delta\lambda_k^{2m} \leq C_2. \label{eig2}
\end{align}

\end{theorem}

\begin{proof}
We show that, for a diagonalizable equation, \eqref{eig1} is equivalent
to \eqref{contA} and \eqref{eig2} is equivalent to \eqref{parab}.
Indeed, note that for every $\gamma,\ r\in \mathbb{R}$,
$$
\|h_k\|_{\gamma+r}=\|\Lambda^rh_k\|_{\gamma}=\lambda_k^r\|h_k\|_{\gamma}.
$$
Then \eqref{eig1} is \eqref{contA} with $v=h_k$, and
\eqref{eig2} is  \eqref{parab} with $v=h_k$. Since both \eqref{eig1} and
\eqref{eig2} are uniform in $k$ and the
collection $\{h_k,\ k\geq 1\}$ is dense in every $\mathbf{H}^{\gamma}$,
 the proof of the theorem is complete.
\end{proof}

The following is the basic existence/uniqueness/regularity result for
 parabolic equations; for the proof,
 see Rozovskii \cite[Theorem 3.2.1]{RozovskiiBook}.

 \begin{theorem}
 \label{th1}
 Assume that equation \eqref{eq2} is parabolic in the triple
 $(\mathbf{H}^{\gamma+m},\mathbf{H}^{\gamma},\mathbf{H}^{\gamma-m})$
 and
 \begin{enumerate}
\item the initial condition $u_0$ is deterministic and belongs to
$\mathbf{H}^{\gamma}$;
\item the process $f=f(t)$ is $\mathcal{F}_t$-adapted with
values in $\mathbf{H}^{\gamma-m}$
and
$$
\mathbb{E}\int_0^T\|f(t)\|_{\gamma-m}^2dt<\infty;
$$
\item each process $g_k=g_k(t)$ is $\mathcal{F}_t$-adapted with values
in $\mathbf{H}^{\gamma}$ and
$$
\sum_{j\geq 1} \mathbb{E} \int_0^T\|g_j(t)\|_{\gamma}^2< \infty.
$$
\end{enumerate}

Then there exists a unique $\mathcal{F}_t$-adapted process $u=u(t)$
with the following properties:
\begin{itemize}
\item $u\in L_2(\Omega; L_2((0,T); \mathbf{H}^{\gamma+m})\bigcap
L_2(\Omega; C((0,T); \mathbf{H}^{\gamma}))$;
\item $u$ is a {\tt solution} of \eqref{eq2}, that is, the equality
\begin{equation*}
\begin{split}
u(t)+\int_0^t(\mathcal{A}_0+\theta\mathcal{A}_1)u(s)ds
&=u_0+\int_0^tf(s)ds\\
&+\sum_{j\geq 1}(\mathcal{M}_ju(s)+g_k(s))dW_j(s).
\end{split}
\end{equation*}
holds in $\mathbf{H}^{\gamma-m}$ for all $t\in [0,T]$ on the
same set $\Omega'\subset\Omega$ of probability one;
\item There exists a positive real number
 $C_0$ depending only on $T$ and the numbers
$C,\delta$ in \eqref{parab} such that
\begin{equation*}
\begin{split}
\mathbb{E}\sup_{0<t<T}\|u(t)\|_{\gamma}^2+
\mathbb{E}\int_0^T\|u(t)\|_{\gamma+m}^2dt&\leq
C_0\left( \|u_0\|_{\gamma}^2+\mathbb{E}\int_0^T\|f(t)\|_{\gamma-m}^2dt
\right.\\
&+
\left.\sum_{j\geq 1} \mathbb{E} \int_0^T\|g_j(t)\|_{\gamma}^2\right)
\end{split}
\end{equation*}
\end{itemize}
\end{theorem}

\begin{corollary}
\label{cor-main}
{\rm
Assume that equation \eqref{eq2} is parabolic and diagonalizable.
Then, under the assumptions of Theorem \ref{th1} we have
\begin{equation}
\label{eq44}
u(t)=\sum_{k=1}^{\infty} u_k(t)h_k {\rm \ and }
\sum_{k=1}^{\infty} \lambda_k^{2\gamma}\mathbb{E}|u_k(t)|^2<\infty,\ t\in [0,T],
\end{equation}
where $u_k(t)=(\Lambda^{\gamma}u(t),h_k)_0$ satisfies
\begin{equation}
\label{eq5}
du_k(t) = \Big((\rho_{k} + \theta\nu_{k})u_k(t) + f_k(t)\Big)dt +
\sum\limits_{j=1}^n (\mu_{jk} u_k(t)+g_k(t)) dW_j(t),
\end{equation}
with $u_k(0) = (\Lambda^{\gamma}u_0,h_k)_0$,
$f_k(t)=\langle \Lambda^{\gamma}f(t),h_k\rangle_{0,m}$,
$g_k(t)=(\Lambda^{\gamma}g(t),h_k)_0$.
}
\end{corollary}

\section{Maximum Likelihood Estimators}\label{section-estimates}

With $(\Omega, \mathcal{F}, \{\mathcal{F}_t\}_{t\geq 0}, \mathbb{P})$,
 $W_j,\ j\geq 1$,  and
$\{\mathbf{H}^r,\ r \in \mathbb{R}\}$ as in the previous section,
consider the stochastic It\^{o} equation
\begin{equation}
\label{eqMain}
du(t)+(\mathcal{A}_0+\theta\mathcal{A}_1)u(t)dt=
\sum_{j\geq 1} \mathcal{M}_ju(t)dW_j(t),\ 0<t\leq T,\ u(0)=u_0.
\end{equation}
We assume that
\begin{itemize}
\item equation \eqref{eqMain} is parabolic in the triple
$(\mathbf{H}^{\gamma+m};\mathbf{H}^{\gamma}, \mathbf{H}^{\gamma-m})$
for some $\gamma\in \mathbb{R}$, $m>0$;
\item equation \eqref{eqMain} is diagonalizable;
\item $u_0\in \mathbf{H}^{\gamma}$.
\item The solution of \eqref{eqMain} is observed (can be measured without
errors) for all $t\in [0,T]$.
\end{itemize}
The objective is to estimate the real number  $\theta$ from the
observations $u(t),\ t\in [0,T]$.

 Even though  whole random field $u$ can be
observed, the actual computations can be performed
 only on a finite-dimensional projection of $u$. By Corollary \ref{cor-main},
 we have
 \begin{align}
 &u(t)=\sum_{k=1}^{\infty}u_k(t)h_k, \label{eq:C1}\\
&u_k(t)+\int_0^t(\rho_k+\theta\nu_k)u_k(s)ds
=(\Lambda^{\gamma}u_0,h_k)_0+\int_0^tu_k(s)\sum_{j\geq 1}\mu_{jk}dW_j(s),
\label{eq:C2}
\end{align}
Thus, a finite collection of the Geometric Brownian motions $u_k$ is a
natural finite-dimen\-sional projection of $u$.

To simplify certain formulas, we will use the following notations:
\begin{equation}
\label{not-eig}
M_k=\sum_{j\geq 1} |\mu_{jk}|^2,\ \eta_k=\frac{M_k}{\nu_k^2}.
\end{equation}

\subsection{Maximum Likelihood Estimator (MLE)}
\label{subsection-MLE}

Let $u_{k_1},\ldots, u_{k_N}$ be a finite collection of diffusion processes
\eqref{eq:C2}. For each $\theta\in \Theta$, the vector
$U_N=(u_{k_1},\ldots, u_{k_N})$ generates a measure on the space of
continuous $\mathbb{R}^N$-valued functions. If these measures are
absolutely continuous with respect to some convenient reference measure, then
the MLE of $\theta$ will be the value maximizing the corresponding
density given the observations. The choice of the reference measure is
dictated, among other factors, by the possibility to find a closed-form
expression of the  density. For diffusion processes with a
parameter in the drift, the standard choice is the measure
generated by the process with  a fixed value of the
parameter, for example, the true value $\theta_0$.  Analysis of the
relevant conditions for mutual absolute continuity, as given,
for example, in the book
by Liptser and Shiryaev \cite[Theorem 7.16]{LiptserShiryayev},
demonstrates that
\begin{itemize}
\item if $N$=1, then the measures generated by $u_k$ for different
values of $\theta$ are mutually
absolutely continuous, and the density with respect to the measure
corresponding to the true parameter $\theta_0$ is
\begin{equation}
\label{eq:dens}
L_k(\theta,\theta_0)=
\exp\Big\{ -\int\limits_{0}^T \frac{\nu_{k}(\theta-\theta_0)}{M_k}
 \, \frac{du_k}{u_k}
-\frac{  \rho_{k} \nu_k  (\theta-\theta_0)\,T}{M_k}
 -
\frac{(\nu_{k})^2(\theta^2-\theta_0^2)\, T\ }{2M_k} \Big\}.
\end{equation}
\item For $N>1$, the measures are typically mutually singular and
so is the resulting estimation problem. We will see later how to exploit this
singularity and gain a computational advantage over the straightforward MLE.
\end{itemize}

Thus, observation of a single process $u_k(t),\ 0\leq t\leq T$,
provides an MLE $\widehat{\theta}_k$ of $\theta$; by
\eqref{eq:dens},
\begin{equation}
\label{eq7}
  \widehat{\theta}_k=  -\frac{1}{\nu_k T} \int\limits_{0}^T\frac{du_k}{u_k}
   - \frac{\rho_k}{\nu_k}.
\end{equation}
By It\^o's Lemma,
$$
d \ln (u_k) = \frac{du_k}{u_k} -\frac12M_kdt,
$$ and
hence from \eqref{eq7} we get
\begin{equation}
\label{eqMLE}
\widehat\theta_k =
\frac{1}{\nu_k T}\ln\frac{u_k(0)}{u_k(T)}
-\frac{M_k}{2\nu_k}.
\end{equation}
Notice that, by uniqueness of solution of equation \eqref{eq:C2},
the function $u_k(t)$ cannot change sign and so $u_k(0)/u_k(T)>0$.
From  \eqref{eq7} and  \eqref{eq:C2} we have the following alternative representation
of the MLE:
\begin{equation}\label{eq8}
  \widehat{\theta}_k = \theta_0 -
  \frac{1}{\nu_kT}\sum_{j\geq1}\mu_{jk}W_j(T);
\end{equation}
in particular,
\begin{equation}
\label{eq:est1}
\mathbb{E}(\widehat{\theta}_k - \theta_0)^2=\frac{\eta_k}{T}
\end{equation}
and
$\sqrt{{T}/{\eta_k}}\,(\widehat{\theta}_k - \theta_0)$
is a standard Gaussian random variable {\em for every $T>0$ and $k\geq 1$.}

All properties of the MLE \eqref{eqMLE} now follow directly from
\eqref{eq8} and \eqref{eq:est1} and are summarized below.

\begin{theorem}
\label{th:MLE}
Assume that equation \eqref{eqMain} is diagonalizable,  parabolic in the triple \\
$(\mathbf{H}^{\gamma+m};\mathbf{H}^{\gamma}, \mathbf{H}^{\gamma-m})$
for some $\gamma\in \mathbb{R}$, $m>0$, and  $u_0\in \mathbf{H}^{\gamma}.$
Then
\begin{enumerate}
 \item For every $k\geq 1$ and $T>0$,
 $\widehat{\theta}_k$ is an unbiased estimator of $\theta_0$.
\item For every $k\geq 1$, as $T\to \infty$,
 $\widehat{\theta}_k$ converges to $\theta_0$ with
 probability one and $\sqrt{T}(\widehat{\theta}_k-\theta_0)$
 converges in distribution to a Gaussian random variable  with zero mean and
 variance $\eta_k$.
\item If, in addition,
\begin{equation}
\label{cconv}
\lim_{k\geq 1} \eta_k=0,
\end{equation}
 then, for every $T>0$,
as $k\to \infty$, $\widehat{\theta}_k$ converges to $\theta_0$ with
 probability one and $(\widehat{\theta}_k-\theta_0)/\sqrt{\eta_k}$
 converges in distribution to a Gaussian random variable  with zero mean and
 variance $1/T.$
\end{enumerate}
\end{theorem}

\begin{remark}{\rm
 Conditions \eqref{eig2} and  \eqref{cconv} are, in general, not
 connected. Indeed, let $\Lambda=\sqrt{I-\boldsymbol{\Delta}}$, where
 $\boldsymbol{\Delta}$ is the Laplace operator on a smooth bounded
 domain in $\mathbb{R}^d$ with zero boundary conditions. Then equation
 $$
 du- (\boldsymbol{\Delta} u  - \theta u)udt=\Lambda udW(t)
 $$
 satisfies \eqref{parab}, but does not satisfy
 \eqref{cconv}: in this case, $\lim_{k\to \infty} \eta_k=\infty$.
Similarly, equation
$$
du-(\theta\boldsymbol{\Delta}u- u)dt= (I-\boldsymbol{\Delta})^{3/4} udW(t)
 $$
 does not satisfy \eqref{parab} for any $\theta$, but satisfies \eqref{cconv}.
 We remark that the solution of this last equation can be constructed in
 special weighted Wiener chaos spaces that are much larger than
 $L_2(\Omega; L_2((0,T); \mathbf{H}^{\gamma}))$; see
 \cite{LR_AP}.
 }
 \end{remark}

 \begin{example}
 {\rm
 Let us consider the following modification of equation \eqref{eq:exmain}
 from Example \ref{ex:main}(b):
 $$
 du- (\boldsymbol{\Delta}u+\theta u)dt=
 \sum_{j\geq 1} (1-\boldsymbol{\Delta})^{-j/2}\, u \,dW_j(t) \ .
 $$
 We have $\nu_k=1$, $\rho_k=-\sigma_k>0$, where $\sigma_k$ are the
 eigenvalues of  $\boldsymbol{\Delta}$, and so $\rho_k\sim ck^{2/d}$;
 $\mu_{jk}=(1+\rho_k)^{-j}$ and
 $$
 M_k=\sum_{j\geq 1} \frac{1}{( 1+\rho_k )^{j}} = \frac{1}{\rho_k} \to 0,\ k\to \infty.
 $$
By Theorem \ref{th:MLE} the maximum likelihood estimator $\widehat{\theta}_k$
 of $\theta$ is
 $$
\widehat{\theta}_k=\frac{1}{T}\ln\frac{u_k(0)}{u_k(T)}+\frac{1}{2\sigma_k}
$$
and
$$
\mathbb{E}(\widehat{\theta}_k-\theta_0)^2\sim cT^{-1}k^{-2/d}
$$
}
\end{example}

\subsection{Modifications of the MLE}
By Theorem \ref{th:MLE}, the MLE \eqref{eqMLE} can be consistent and
asymptotically normal either in the limit $T\to \infty$ or in the
limit $k\to \infty$. An increase of $T$ always improves the
quality of the estimator by reducing the variance; if
\eqref{cconv} holds, then the variance of the estimator can be
further reduced by using $u_k$ with the largest available value
$k$.

The natural question is whether the quality of the
estimator can be improved even more by using more than one
process $u_k$. This question is no longer of statistical nature:
as equation \eqref{eq:C2} shows, each $u_k$ contains essentially  the
same stochastic information. More precisely,   the sigma-algebra
generated by each $u_k(t)$, $t\in [0,T]$ coincides with the
sigma-algebra generated by $\mu_{jk}W_j(t),\ j\geq 1,\ t\in [0,T]$
(some of $\mu_{jk}$ can, in principle, be zeroes).
Moreover, as was mentioned above, the statistical estimation model
for $\theta$, involving two or more processes $u_k$, is singular.
In what follows,  we will see how to use this singularity to gain computational
advantage over \eqref{eqMLE}.

The problem can now be stated as follows: given a sequence of numbers
$\widehat{\theta}_k$ such that $\lim_{k\to \infty} \widehat{\theta}_k=\theta_0$,
can we transform it into a sequence $\widetilde{\theta}_k$  such  that
\begin{equation}
\label{bconv}
\lim_{k\to \infty} \widetilde{\theta}_k=\theta_0,\
\limsup_{k\to \infty}
\frac{|\widetilde{\theta}_k-\theta_0|}{|\widehat{\theta}_k-\theta_0|}<1.
\end{equation}
If \eqref{bconv} holds, it is natural to say that
$\widetilde{\theta}_k$ converges to $\theta_0$ faster than
$\widehat{\theta}_k$. Accelerating the convergence
of a sequence is a classical problem in numerical analysis. The
main features of this problem are (a) There are many different
methods to accelerate the convergence, and (b) the effectiveness of every
method varies from sequence to sequence.

We will investigate two methods:
\begin{enumerate}
\item Weighted averaging;
\item Aitken's $\triangle^2$ method.
\end{enumerate}

\begin{theorem}[Weighted averaging]
Let $\beta_k,\ k\geq 1,$ be a sequence of non-negative numbers and
$$
\sum_{k\geq1}\beta_k=+\infty.
$$
Define the weighted averaging estimator $\widehat\theta_{(N)}$ by
\begin{equation}
\label{eq10}
\widehat\theta_{(N)} = \frac{\sum\limits_{k=1}^{N}
\beta_k\widehat\theta_k }{\sum\limits_{k=1}^{N}\beta_k} \, .
\end{equation}
Then
\begin{enumerate}
 \item For every $N\geq 1$ and $T>0$,
 $\widehat{\theta}_{(N)}$ is an unbiased estimator of $\theta_0$.
\item For every $N\geq 1$, as $T\to \infty$,
 $\widehat{\theta}_{(N)}$ converges to $\theta_0$ with
 probability one and $\sqrt{T}(\widehat{\theta}_{(N)}-\theta_0)$
 converges in distribution to a Gaussian random variable  with zero mean and
 variance
 \begin{equation}
 \label{WVar1}
V_N= \sum_{j\geq 1} \left(\frac{\sum_{k=1}^N(\beta_k \mu_{jk}/\nu_k)}
 {\sum_{k=1}^N\beta_k}\right)^2.
 \end{equation}
\item If, in addition, \eqref{cconv} holds
 then, for every $T>0$,
as $N\to \infty$, $\widehat{\theta}_{(N)}$ converges to $\theta_0$ with
 probability one.
\end{enumerate}
\end{theorem}

\begin{proof}
By  \eqref{eq8},
\begin{equation}
\label{eq11}
\widehat\theta_{(N)}  =
\theta_0 +
\frac{ \sum\limits_{j\geq 1}\left(
 \sum\limits_{k=1}^{N} (\beta_k \mu_{jk}/\nu_k)\right) W_j(T)}
{T \sum\limits_{k=1}^N\beta_k},
\end{equation}
from which  the first two statement of the theorem follow.
 For the last statement, we combine \eqref{eq10} with
  the Toeplitz lemma: if $\lim_{k\to \infty} a_k=a$ and $\beta_k>0$,
then
$$
\lim_{N\to \infty} \frac{\sum_{k=1}^N \beta_ka_k}{\sum_{k=1}^N\beta_k}=a.
$$
\end{proof}

 The behavior of $V_N/\eta_N$, as $N\to \infty$ can be just about
 anything.
 Take  $\rho_k=0$,  $\mu_{jk}=0$, $j>1$. Then,
 \begin{itemize}
\item With $\beta_k=1/k$, $\nu_k=k^2$, and $\mu_{1k}=k,$
  we get  $\eta_N=1/N^2$ and
$$
\frac{V_N}{\eta_N}\sim \frac{\alpha N^2}{\ln^2N}\to \infty, \ N\to \infty.
$$
for some $\alpha>0$; recall that,  for $a_n,b_n>0$, notation
$a_n\sim b_n$ means
$$
\lim_{n\to \infty} (a_n/b_n)=1.
$$
\item  With $\beta_k=k$, $\nu_k=k^2$, and $\mu_{1k}=k,$
  we get  $\eta_N=1/N^2$ and
$$
V_N\sim 4\eta_N>\eta_N
$$
\item With $\beta_k=1$, $\nu_k=k^2$, and $\mu_{1k}=(-1)^kk,$
  we get  $\eta_N=1/N^2$ and
  $$
  V_N\sim (\ln^22)\eta_N<\eta_N
$$
\item With $\beta_k=1$, $\nu_k=k$, and $\mu_{1k}=(-1)^k\sqrt{k},$
  we get  $\eta_N=1/N$ and
$$
\frac{V_N}{\eta_N}\sim \frac{\beta}{{N}}\to 0,\ N\to \infty.
$$
\end{itemize}

Next, we consider {\bf  Aitken's $\triangle^2$ method}. This method
consists in transforming a sequence
$A=\{a_n,\ n\geq 1\}$ to a sequence
$$
b_n(A)=a_n-\frac{(a_{n+1}-a_n)^2}{a_{n+2}-2a_{n+1}+a_{n}}.
$$
The main result concerning this method is that if
$\lim_{n\to \infty} a_n=a$ and
\begin{equation}
\label{A1}
\lim_{n\to \infty}\frac{|a_{n+1}-a|}{|a_n-a|}=\lambda \in (0,1),
\end{equation}
then $\lim_{n\to \infty}b_n(A)=a$ and
$$
\lim_{n\to \infty} \frac{|b_n(A)-a|}{|a_n-a|}=0.
$$
That is, the sequence $b_n(A)$ converges to the same limit $a$ but faster.

Accordingly, under the condition \eqref{cconv}, we define
\begin{equation}
\widetilde{\theta}_k=
\widehat{\theta}_k-\frac{(\widehat{\theta}_{k+1}-\widehat{\theta}_k)^2}{
\widehat{\theta}_{k+2}+2\widehat{\theta}_{k+1}-\widehat{\theta}_k},
\end{equation}
with a hope that
\begin{equation}
\label{A2}
\lim_{k\to \infty}\frac{\mathbb{E}(\widetilde{\theta}_k-\theta_0)^2}
{\mathbb{E}(\widehat{\theta}_k-\theta_0)^2}<1.
\end{equation}
In general, there is no guarantee that this will be the case
because typically $\eta_k\sim \alpha k^{-\delta}$ for some $\alpha>0$ and
 $\delta>0$, and so, if we set
$$
a_k=\mathbb{E}(\widehat{\theta}_k-\theta_0)^2,
$$
we get by Theorem \ref{th:MLE}
$$
\lim_{n\to \infty}\frac{|a_{n+1}-a|}{|a_n-a|}=1.
$$
Direct investigation of the sequence $\widetilde{\theta}_k$ is possible
if there is only one Wiener process $W=W(t)$ driving the equation, that
is, $\mu_{jk}=0$ for $j\geq 2, \ k\geq 1.$ In this case, \eqref{eq8} shows that
\begin{equation}
\label{A3}
\widetilde{\theta}_k=\theta_0+\frac{W(T)}{T}\,\left(
r_k-\frac{(r_{k+1}-r_k)^2}{r_{k+2}-2r_{k+1}+r_{k}}\right),
\end{equation}
where $r_k=\mu_{1k}/\nu_k$.
Then direct computations show that
\begin{itemize}
\item if $r_k\sim\alpha k^{-\delta}$, $\alpha, \, \delta>0$, then
$$
\frac{\mathbb{E}(\widetilde{\theta}_k-\theta_0)^2}
{\mathbb{E}(\widehat{\theta}_k-\theta_0)^2} \sim \frac{1}{(\delta+1)^2}.
$$
\item if $r_k=(-1)^k/k$, then
$$
\frac{\mathbb{E}(\widetilde{\theta}_k-\theta_0)^2}
{\mathbb{E}(\widehat{\theta}_k-\theta_0)^2} \sim \frac{c}{k^2},\ c>0.
$$
\end{itemize}

For more than one Wiener process, we find
$$
\widetilde{\theta}_k=\theta_0+\frac{\xi^2_k}{\zeta_k},
$$
where $(\xi_k, \zeta_k)$ is a two-dimensional  Gaussian vector with
known distribution. The analysis of this estimator, while possible, is
technically much more difficult and will require many additional assumptions
on $\mu_{jk}$. We believe that this analysis falls outside the scope of this paper,
and we present here only some numerical results. We suppose that Fourier coefficients
$u_k$ satisfy \eqref{eq:C2}   with
$\nu_k=k, \rho_k=0, \ \mu_{jk}=(-1)^k/(k+j)$, the noise term is driven  by
$n=10$ Wiener processes, and the true value of the parameter $\theta_0=1$.
From \eqref{eqMLE} we note that the estimates $\hat{\theta}_k$
 can be calculated if we only know the value of $\log(u_k(T)/u_k(0))$, rather than the whole path
$u_k(t), \ 0\leq t\leq T$.
Using the closed-form solution of equation \eqref{eq:C2}
$u_k(t) = u_k(0)\exp(-(\theta_0\nu_k + \sum_{j}\mu_{jk}^2/2)t + \sum_{j}\mu_{jk}W_{j}(t))$,
we simulate $\log(u_k(T)/u_k(0))$ directly, without applying some discretization schemes to the process $u_k(t)$.
Three type of estimates are presented in Figure 1. The obtained numerical results are consistent with above
theoretical results: Aitken's $\triangle^2$ method performs the best, Weighted Averages Estimates
with $\beta_k=k$ perform better than simple estimates.

\begin{figure}[ht]
\centering
\includegraphics[width=.95\linewidth]{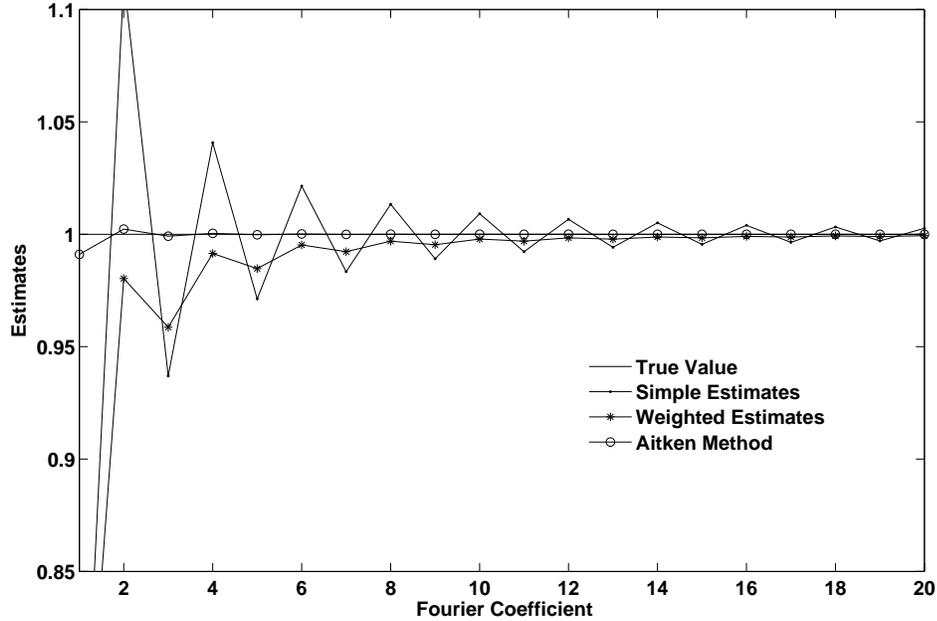}
\caption{ {\footnotesize Performances of three type of estimates: Simple, Weighted Averages
and Aitken's $\triangle^2$ method} }
\end{figure}

\section{Closed-form Exact  Estimators}\label{subsection-exact}
In regular models, the estimator is consistent in the large sample or
small noise limit; neither of these limits can be evaluated exactly  from any actual
observations.
In singular models, there often exists an estimator that is consistent
in the limit that can potentially be evaluated exactly  from the available
observations. Still, no expression can be evaluated on
a computer unless the expression involves only finitely many operations of
addition, subtraction, multiplication, and division.

\begin{definition} An estimator is called {\tt closed-form exact} if
it produces the exact value of the unknown parameter after a finite number
of additions, subtractions, multiplications, and divisions performed on the
elementary functions of the observations.
\end{definition}

Closed-form exact estimators  exist for
the model \eqref{eqMain} {\em if we assume that the observations are
$u_k(t),\ k\geq 1,\ t\in [0,T]$.}

As an illustration, consider the simple example
$$
du-\theta u_{xx}dt = (u/2)dt+ udW(t),
$$
where $x\in (0,\pi)$ and zero boundary conditions are assumed.

With $h_k=\sqrt{2/\pi}\,\sin(kx)$, we find
$$
du_k(t)=-k^2\theta u_k(t)dt+(u_k/2)dt+u_k(t)dW(t).
$$
Set $v_k(t)=\ln (u_k(t)/u_k(0))$. Then
$$
dv_k(t)=-k^2\theta dt+dW(t).
$$
In particular,
$$
v_1(T)=-\theta T+W(T),\ v_2(T)=-4\theta T + W(t)
$$
so that
$$
\theta=\frac{v_1(T)-v_2(T)}{3T}
$$
or
\begin{equation}
\label{EXP1}
\theta=\frac{1}{3T}\ln\frac{u_1(T)u_2(0)}{u_1(0)u_2(T)}.
\end{equation}
Notice that given $u_1, \ldots, u_N$, we have $N(N-1)/2$ exact estimators
of this type.

If there are two Wiener processes driving the equation, then we will
need three different $u_k$ to construct an estimator of the type
\eqref{EXP1}. The general result is as follows.

\begin{theorem}
In addition to conditions of Theorem \ref{th:MLE} assume that
there exist two finite sets of indices
$(k_1^i, k_2^i,\ldots, k_n^i), \ i=1,2$, and a positive integer $p$ such that
 $$
 \sum_{\ell= 1}^n\nu_{k_{\ell}^1}\mu_{pk_{\ell}^2}\neq
\sum_{\ell= 1}^n\nu_{k_{\ell}^2}\mu_{pk_{\ell}^1}.
$$
Then there exists a closed-form exact estimator for $\theta$.
\end{theorem}

\begin{proof}
Let $v_k(t)=\ln (u_k(t)/u_k(0))$. From \eqref{eq:C2}, by It\^o's formula, we get
\begin{equation}\label{eq9}
d v_k = -\Big(\rho_k + \theta\nu_k + \frac{1}{2}M_k\Big)dt +
\sum\limits_{j\geq1}\mu_{jk}dW_j(t),
\end{equation}
and then
\begin{equation}
\label{eq9_3}
\begin{split}
\theta\sum_{\ell= 1}^n\Big(\nu_{k_{\ell}^2}\mu_{pk_{\ell}^1}-
\nu_{k_{\ell}^1}\mu_{pk_{\ell}^2}\Big)&=
\frac{1}{T}\sum_{\ell= 1}^n\Big(\mu_{pk_{\ell}^2}
v_{k_{\ell}^1}(T)-\mu_{pk_{\ell}^1}
v_{k_{\ell}^2}(T)
\\ &+
\mu_{pk_{\ell}^2}(\rho_{k_{\ell}^1}+\frac{1}{2}M_{k_{\ell}^1})-
\mu_{pk_{\ell}^1}(\rho_{k_{\ell}^2}+\frac{1}{2}M_{k_{\ell}^2})
\Big),
\end{split}
\end{equation}
which completes the proof.
\end{proof}

If there are $n_0$ Wiener processes driving the equation, then
the extra condition of the theorem can always be ensured with
$n=n_0+1$, because every collection of $n$ vectors in an
$n-1$-dimensional space is linearly dependent.
While relation \eqref{eq9_3} gives a closed-form exact estimator, the
resulting formulas can be rather complicated when the number
of Wiener processes in the equation is large; if this number is infinite, then
the estimator might not exist at all. For comparison, the complexity of the
 maximum likelihood estimator \eqref{eqMLE} does not depend on the number
of Wiener processes in the equation. As a result, when it comes to
actual computations, the closed-form exact
estimator is not necessarily the best choice. On the other hand, the
very existence of such an estimator is rather remarkable.

We conclude this section with three examples of
 closed-form exact estimators.
The first example shows that such estimators can exist for equations
 that are not diagonalizable in the sense of Definition \ref{def000}.

 \begin{example}
 {\rm
 Consider the equation
 $$
 du(t,x)= \theta u_{xx}(t,x)dt +u(t,x)dW(t),\ 0<t\leq T,\ x\in \mathbb{R}.
 $$
 By the It\^{o} formula,
 $$
 u(t,x)=v(t,x)\exp(W(t)-(t/2)),
 $$
  where $v$ solves the heat equation
 $v_t=\theta v_{xx}$, $v(0,x)=u(0,x)$.
Assume that $u(0,x)$ is a smooth compactly supported function.
Then $u(t,x)$ is a smooth bounded function for all $t>0, x\in \mathbb{R}$
and $\mathbb{E}\int_{\mathbb{R}}|u(t,x)|^pdx<\infty$ for all $p>0$, $t\geq 0$.
In particular, the Fourier transform $U(t,y)$  of $u$ is defined and
satisfies
$$
dU(t,y)=-\theta y^2U(t,y)dt+U(t,y)dW(t).
$$
Let $V(t)=\ln (U(t)/U(0))$. Then
$$
V(T,y)= -y^2\theta T-(T/2)+W(T),
$$
and
$$
\theta=\frac{V(T,y_1)-V(T,y_2)}{T(y_2^2-y_1^2)}.
$$
}
\end{example}

The next example shows that conditions \eqref{parab} and \eqref{cconv}
 are not related  to the
existence of a closed-form exact estimator.

\begin{example}
{\rm
Consider the equation
$$
du-(\boldsymbol{\Delta}u+\theta u)dt= (I-\boldsymbol{\Delta})^{3/4} udW(t)
 $$
 on $(0,\pi)$ with zero boundary conditions. Clearly both
 \eqref{parab} and \eqref{cconv} are not satisfied.
While the equation is not parabolic, there exists a unique solution in
weighted Wiener chaos spaces, and we can therefore consider
$$
du_k=(-k^2u_k+\theta u_k)dt-(1+k^2)^{3/4}u_kdW(t).
$$
For $v_k(t)=\ln (u_k(t)/u_k(0))$ we find
$$
v_k(T)=(-k^2-\frac{(1+k^2)^{3/2}}{2})T+\theta T + (1+k^2)^{3/4}W(t).
$$
In particular,
$$
v_1(T)=a_1T+b_1W(t)+\theta T,\ v_2(T)=a_2T+b_2W(T)+\theta T,
$$
and so
$$
\theta=\frac{b_1v_2(T)-b_2v_1(T)-(a_2b_1-a_1b_2)T}{T(b_1-b_2)}.
$$
}
\end{example}

The last example shows that, as long as there is no spacial structure in the
noise,  multiplicativity  of the noise is not
necessary to have a closed-form exact estimator.
\begin{example}
{\rm
Consider the equation
$$
du(t,x)=\theta u_{xx}(t,x)dt + dW(t), t>0,\ x\in (0,\pi),
$$
with Neumann boundary conditions,
so that $h_1=1/\sqrt{\pi}$ and $h_k=\sqrt{2/\pi}\cos((k-1)x)$,
$k\geq 2$.
Then $du_2(t)=-\theta u_2(t)dt$, and, as long as
$u_2(0)\not=0$, we have
$$
\theta=\frac{1}{T}\ln\frac{u_2(0)}{u_2(T)}.
$$
}
\end{example}

\def\cprime{$'$} \def\cprime{$'$} \def\cprime{$'$}
\providecommand{\bysame}{\leavevmode\hbox to3em{\hrulefill}\thinspace}
\providecommand{\MR}{\relax\ifhmode\unskip\space\fi MR }
\providecommand{\MRhref}[2]{%
  \href{http://www.ams.org/mathscinet-getitem?mr=#1}{#2}
}
\providecommand{\href}[2]{#2}

\end{document}